\documentclass[reqno]{amsart}
\usepackage{amsmath}
\usepackage{amssymb}
\usepackage{amsthm}
\usepackage{enumerate}
\usepackage[mathscr]{eucal}
\usepackage{graphicx}

\newtheorem{theorem}{Theorem}[section]
\newtheorem{lemma}[theorem]{Lemma}
\newtheorem{proposition}[theorem]{Proposition}
\newtheorem{corollary}[theorem]{Corollary}
\newtheorem{J-S Theorem}[theorem]{J-S Theorem}
\newtheorem{definition}[theorem]{Definition}
\theoremstyle{definition}
\newtheorem{example}[theorem]{Example}
\newtheorem{remark}[theorem]{Remark}

\begin{document}

		\title{Left K-Cauchy regular maps between quasi-metric spaces}
	\author{Om Dev Singh and Anubha Jindal \textsuperscript {}}
	
	\thanks{Om Dev Singh: Department of Mathematics, Malaviya National Institute of Technology Jaipur, Jaipur, Rajasthan, India, email: om2947dev@gmail.com}
	\thanks{Anubha Jindal: Department of Mathematics, Malaviya National Institute of Technology Jaipur, Jaipur, Rajasthan, India, email: anubha.maths@mnit.ac.in}
		
		
		\maketitle
			\begin{abstract}
			This article presents a systematic study of a class of maps between quasi-metric spaces that preserve left K-Cauchy sequences. We call such maps left K-Cauchy regular maps. Several characterizations of these maps have been given in terms of the structural properties of the underlying quasi-metric spaces. In particular, we characterize left K-Cauchy regular maps in terms of hereditarily precompact sets by using a suitable version of the classical Efremovi$\acute{c}$ Lemma in the setting of quasi-metric spaces. In addition, we examine the relationship of left K-Cauchy  regular maps with (forward) uniformly continuous and continuous maps and prove extension theorems for such maps. 		
		\end{abstract}

	\section{Introduction}
	A quasi-metric $d$ on a set $X$ is a function on $X\times X$ that satisfies all the metric axioms, but the symmetry condition $d(x,y) =d(y,x)$ may fail for $x,y \in X$. This subtle yet significant modification allows quasi-metric spaces to more accurately model a variety of real-world phenomena where direction plays a crucial role, such as in transportation networks, communication systems, and biological processes. These spaces are also sometimes called asymmetric metric spaces or oriented metric spaces \cite{Bodjanova,Zimmer}.
		
	For a quasi-metric space $(X,d)$, there is another naturally associated quasi-metric $\bar{d}$ on $X$ defined as $\overline{d}(x,y) = d(y,x)$ for $x, y \in X$. This $\bar{d}$ is called the conjugate quasi-metric of $d$. Due to the presence of this dual structure, many topological and metric concepts, such as completeness, compactness, convergence, continuity, and uniform continuity for quasi-metric spaces, extend beyond the usual analysis of metric spaces. The monograph by  Fletcher and Lindgren \cite{fletcherlindgren} presents the fundamental theory of quasi-uniform and quasi-metric spaces. Cobza\c s \cite{Cobzas} presents a modern approach to quasi-metric spaces. In this monograph, Cobza\c s treated the theory of asymmetric normed spaces in detail and gave asymmetric counterparts of several important theorems of classical functional analysis, such as Hahn-Banach and Krein-Milman type theorems and classical open mapping and closed graph theorems.

	In recent years, several researchers have shown interest in developing analysis on quasi-metric spaces, especially with respect to function spaces, hyperspaces, and fixed point theory \cite{RomagueraPedro_fp, cow1, kunjiromaguera, RomagueraLopez,Seceleanetal_fp}. For more on quasi-uniform spaces and quasi-metric spaces, we refer readers to survey articles by K\"unzi \cite{kunji3,kunji4,kunji2}.

	Despite the above advances, the study of various classes of functions between quasi-metric spaces remains largely unexplored. It is well-known that several structural properties of metric spaces can be characterized in terms of functions defined on them. For instance, a metric space is compact if and only if every real-valued continuous function on it is bounded. Consequently, various classes of functions between metric spaces have been defined and studied in detail \cite{beergariddo2020,Beer_Stability,borsik1,Lipsy,snipes1}.
	
	One important class of such functions is the class of Cauchy regular functions. A function between two metric spaces is said to be Cauchy regular if it preserves Cauchy sequences.  In \cite{snipes1,snipes3}, Snipes defined and examined this class, which lies strictly between the classes of uniformly continuous functions and continuous functions. He proved that a metric space is complete if and only if every real-valued continuous function defined on it is Cauchy regular. 	Bors\'ik further explored these functions in \cite{borsik1, borsik2}. Since every Cauchy regular function is continuous, these functions are sometimes referred to as Cauchy continuous functions in the literature \cite{aggarwalkundu, beergarrido2,beergariddo2020,Beer_totalBounded, kundumanisha}. In \cite{beergarrido2014}, Beer and Garrido proved that a metric space is totally bounded if and only if each real-valued Cauchy continuous map defined on it is bounded, and in \cite{beergarrido2}, they proved that the class of Cauchy-Lipschitz functions is uniformly dense in the class of Cauchy continuous functions. For further details on this topic, we refer to recent research monographs by Beer \cite{beer1} and Kundu et al. \cite{kundulipsymanisha}.

	In this paper, we present a systematic study of Cauchy regular functions between quasi-metric spaces, aiming to develop their foundational theory and structural properties within the broader framework of these asymmetric spaces. Here we would like to emphasize that the study of Cauchy regular functions between quasi-metric spaces is far from straightforward. The absence of symmetry not only alters several classical properties but also demands new techniques and insights to obtain analogous results available for metric spaces.
	
	In the literature on quasi-metric spaces, several non-equivalent notions of Cauchy sequences have been investigated in order to study the completeness of these spaces \cite{Cobzas_Ekeland, doitchinov1988,kelley1963,reillysubramurthy1982, roma_leftcompletion}. Most of these classes of Cauchy sequences coincide in the realm of metric spaces. Among the various notions of Cauchy sequences, the most studied notion is that of a left K-Cauchy sequence.  A sequence $(x_n)$ in a quasi-metric space $(X,d)$ is called \textit{left K-Cauchy} (\textit{right K-Cauchy}) provided for each $\epsilon >0$ there exists $n_0 \in \mathbb{N}$ such that for $n\geq k \geq n_0$, we have $d(x_k, x_n) < \epsilon$ ($d(x_n, x_k) < \epsilon$). A function between two quasi-metric spaces that preserves left K-Cauchy sequences is called \text{left K-Cauchy regular}. The behavior of left K-Cauchy regular maps between quasi-metric spaces significantly differs from that of Cauchy regular maps between metric spaces. For instance, every Cauchy regular map is continuous, but it need not be true for a left K-Cauchy regular map (see Example \ref{e10}).

	The organization of this paper is as follows. In Section 2, we discuss some basic definitions and results related to quasi-metric spaces. Section 3 presents a systematic study of left K-Cauchy regular maps between quasi-metric spaces. In particular, we provide a sequential characterization of these maps. Then, by using a suitable version of the classical Efremovi$\acute{c}$ Lemma in the setting of quasi-metric spaces, we characterize left K-Cauchy regular maps in terms of hereditarily precompact sets (see, Theorem \ref{distance_leftKCauchy}). As a corollary, we obtained that if a map $f: (X,d) \to (Y,\rho)$ between two quasi-metric spaces $(X,d)$ and $(Y,\rho)$ is left K-Cauchy regular, then it is Cauchy regular as a map between the associated metric spaces. However, the converse fails.  In the last two sections, we examine the relation of left K-Cauchy regular maps with continuous and uniformly continuous functions and prove extension theorems for such maps.  Under certain assumptions, it is shown that a map between two quasi-metric spaces is left K-Cauchy regular if and only if its composition with every real-valued uniformly continuous function is left K-Cauchy regular. Several examples and counterexamples are given to support our results.

	\section{Preliminaries}
	In this section, we summarize some basic definitions and results related to quasi-metric spaces. For more details readers may refer to \cite{Cobzas, Zimmer}.
	
	\begin{definition}\normalfont
		A \textit{quasi-metric} on a set $X$ is a function $d: X \times X \rightarrow \mathbb{R}$ satisfying the following properties for $x, y, z \in X$
		
		\begin{enumerate}
			\item  $d(x, y) \geq 0$, $d(x,x)=0$
			\item  $d(x, y) = 0$ = $d(y,x)$ $\implies$ $x = y$	
			\item  $d(x, y) \leq d(x, z)+ d(z, y)$.
		\end{enumerate}
	\end{definition}
	
	If $d$ is a quasi-metric, then the \textit{conjugate quasi-metric} $\bar{d}$ is defined as $\overline{d}(x,y)=d (y,x)$ for $ x,y\in \;X.$
	Moreover, $d^{s}(x,y)=\max\{d(x,y),\overline{d}(x,y)\}$ for $x, y \in X$ is a metric on $X$, this metric is called the \textit{associated metric} to the quasi-metric $d$.
	
	\begin{definition}\normalfont
		Suppose $(X,d)$ is a quasi-metric space. Then
		\begin{enumerate}
			\item $B^{+}(x,\epsilon)=\{y\in X:d(x,y)<\epsilon\}$ is called the \textit{forward open ball}; and
			\item $B^{-}(x,\epsilon)=\{y\in X:d(y,x)<\epsilon\}$ is called the \textit{backward open ball}.
		\end{enumerate}  
	\end{definition}
	The topology $\tau(d)$ (resp. $\tau(\overline{d})$) on $X$ generated by all forward balls (resp. backward balls) is called the \textit{forward topology} (resp. \textit{backward topology}) induced by $d$. Clearly, every quasi-metric space is first countable.

	\begin{example}\label{e4}
		The function $d:\mathbb{R}\times \mathbb{R} \to \mathbb{R}$  defined by 
		$d(x,y)=\begin{cases}
			y-x, \text{ if } y\geq x\\    
			0, \text{   if }  y<x
		\end{cases}$\\
		is a quasi-metric on $\mathbb{R}$. Here $B^{+}(x,\epsilon)=(-\infty,x+\epsilon)$ and $B^{-}(x,\epsilon)= (x-\epsilon,\infty)$.
		The quasi-metric $d$ is called the \textit{upper quasi-metric} and its conjugate is called the \textit{lower quasi-metric} on $\mathbb{R}$. The associated metric is the usual metric on $\mathbb{R}$. \qed 
	\end{example}
	The space of reals with the usual  metric is denoted by $(\mathbb{R},|\cdot|)$.
	\begin{example}\label{e3} 
		The function	$d :\mathbb{R} \times \mathbb{R} \rightarrow  \mathbb{R}$ defined by
		$ d(x,y)=\begin{cases}
			y-x, \text{ if } y\geq x\\    
			1, \text{   if }  y<x
		\end{cases}$\\
		is a quasi-metric on $\mathbb{R}$. The forward topology on $\mathbb{R}$ induced by $d$ is the well-known lower limit or Sorgenfrey topology. Consequently, we call this quasi-metric  on $\mathbb{R}$, the \textit{Sorgenfrey quasi-metric}. \qed
	\end{example}

	\begin{remark}The following facts show how the theory of quasi-metric spaces is different from metric spaces. 
		\begin{enumerate}

			\item It is easy to see that every quasi-metric space is $T_0$, however they need not be $T_1$ (see, Example \ref{e4}).
						\item A quasi-metric space $(X,d)$ is $T_1$ if and only if $d(x,y)>0,$ whenever $x\neq y$. Also, $\tau(d)$ is $T_1$ if and only if $\tau(\bar{d})$ is $T_1$ (\cite{Cobzas}).
			
			\item  From Example \ref{e3}, it follows that separability and second countability are not necessarily equivalent in quasi-metric spaces.			
		\end{enumerate}
	\end{remark}
	In the sequel, several facts are given that further demonstrate the difference between the theory of quasi-metric and metric spaces.  	
	
	\begin{definition}\normalfont
		A sequence $(x_{n})$ in a quasi-metric space $(X,d)$ is said to be \textit{forward convergent} (resp. \textit{backward convergent}) if there exists $ x\in X$ such that $\lim_{n}d(x,x_{n})=0$ (resp. $\lim_{n}d(x_{n},x)=0$). If $(x_n)$ is forward convergent (resp. backward convergent) to $x$, we denote it by $x_{n}\xrightarrow{f}x$ (resp. $x_{n}\xrightarrow{b}x$).
	\end{definition}
	\begin{remark}\label{R1} 
		In quasi-metric spaces forward convergence may not imply backward convergence. For instance in the Sorgenfrey quasi-metric space (Example \ref{e3}), the sequence $(1/n)$ is forward convergent to $0$ but it is not backward convergent to $0$.  In a $T_1$ quasi metric space $(X,d)$, if a sequence $(x_n)$ is forward convergent to $x\in X$ and backward convergent to $y\in X$, then $x=y$. In particular, if in a $T_1$ quasi metric space $(X,d)$, forward convergence implies backward convergence, then the forward limit is unique (see, Lemma 3.1 in \cite{Zimmer}). 
	\end{remark}
	
	Quasi-metric spaces having the property forward convergence implies backward convergence are also referred to as small set symmetric in literature \cite{kunzireilly}. 
	
	\begin{example}
		Let \( d : \mathbb{R} \times \mathbb{R} \rightarrow \mathbb{R}^{+}\cup\{0\} \) be defined by
		\[
		d(x, y) = 
		\begin{cases} 
			d(x, 0) + d(0, y), & \text{if } y \neq x \\
			0, & \text{if } x = y,
		\end{cases}
		\]
		where \( d(x, 0) = \frac{1}{|x|} \) and \( d(0, x) = \frac{1}{|x|^2} \) for \( x \neq 0 \). Then $d$ is a quasi-metric on $\mathbb{R}$.
		
		Let \( 0 \neq x \in \mathbb{R} \). If a sequence $(x_n)$ is forward convergent to $x$, that is, $ x_n \xrightarrow{f} x$, then there exists \( N \in \mathbb{N} \) such that \( x_n = x \) for \( n \geq N \). Therefore, \( x_n \xrightarrow{b} x \). Similarly, if \( x = 0 \) and \( x_n \xrightarrow{f} 0 \), then \( |x_n| \rightarrow \infty \) in $(\mathbb{R},|\cdot|)$. Therefore,  \( x_n \xrightarrow{b} 0 \). Thus, the quasi-metric space $(X,d)$ has the property forward convergence implies backward convergence. \qed
	\end{example}

	\begin{definition}\normalfont Let $(X,d)$ be a quasi-metric space. A sequence $(x_n)$ in $(X,d)$ is called \textit{left K-Cauchy} (resp. \textit{right K-Cauchy}) if for every $\epsilon>0$ there exists  $n_{0} \in \mathbb{N}$ such that $\forall$ $ n_{0} \leq k \leq n \implies d(x_{k},x_{n})<\epsilon$ (resp. $d(x_{n},x_{k})<\epsilon$).
	\end{definition}
	
	\begin{remark}\label{conv_not Cauchy}
		In metric spaces every convergent sequence is Cauchy, but in a quasi-metric space a forward convergent sequence need not be left K-Cauchy. For instance, in the Sorgenfrey quasi-metric space (Example \ref{e3}), the sequence $(1/n)$ is forward convergent to $0$ but it is not left K-Cauchy. However, it can be easily proved that if a left K-Cauchy sequence $(x_n)$ has a forward convergent subsequence, then $(x_n)$ will also be forward convergent.
		
	\end{remark}
	\begin{definition}\normalfont
		A quasi-metric space $(X,d)$ is called \textit{left K-complete} if every left K-Cauchy sequence is forward convergent.
	\end{definition}
	
	\begin{example}
		Consider the quasi-metric space $(\mathbb{R},d)$, where $d$ is the upper quasi-metric as defined in Example \ref{e4}. Then $(\mathbb{R},d)$ and $(\mathbb{R},\bar{d})$  are left K-Complete. However, the Sorgenfrey quasi-metric space (Example \ref{e3}) is not left K-complete as $(x_n)=(\frac{-1}{n})$ is left K-Cauchy but is not forward convergent. \qed
	\end{example}
		
		\begin{definition}\normalfont	A function $f:(X,d)\longrightarrow (Y,\rho)$ is called forward-forward continuous (ff-continuous), if for every $\epsilon >0$, there exists $\delta >0$ such that $y\in B^+(x,\delta)$ implies $f(y)\in B^+(f(x),\epsilon)$.
		\end{definition}
		Similarly, one can define fb-continuous, {bf-continuous} and {bb-continuous}.
				\begin{definition}\normalfont Let $(X,d)$ and $(Y,\rho)$ be quasi-metric spaces.
				A function $f:(X,d)\rightarrow(Y,\rho)$ is called \textit{forward uniformly continuous} (resp. \textit{backward uniformly continuous}) if for every $\epsilon >0,$ there exists $\delta>0$ such that $\rho(f(x),f(y))<\epsilon$ (resp. $\rho(f(y),f(x))<\epsilon$) whenever $d(x,y)<\delta$.
			\end{definition}
			Clearly, every forward uniformly continuous function is both ff-continuous and bb-continuous; and every backward uniformly continuous function is both fb-continuous and bf-continuous. Also it is easy to see that a forward uniformly continuous map between quasi-metric spaces  maps a left K-Cauchy sequence to a left K-Cauchy sequence.
			
						\begin{definition}\normalfont (Definition 3.1, \cite{Moshokoa}) Let $(x_n)$ and $(y_n)$ be sequences in a quasi-metric space $(X,d)$. Then $(x_n)$ is said to be \textit{forward parallel} (resp. \textit{backward parallel}) to $(y_n)$ denoted as $ (x_n){\overset{f}{||}} (y_n)$ (resp. $ (x_n)\overset{b}{||} (y_n)$) if for every $ \epsilon >0$, there exists
				$ N \in \mathbb{N}$ such that $d(x_n,y_n) <\epsilon$ (resp. $d(y_n,x_n) <\epsilon$) whenever $ n\geq N$.

				If $(x_n)$ is both forward and backward parallel to $(y_n)$, then we denote it by $(x_n)\overset{f}{\underset{b}{||}} (y_n)$.
			\end{definition}
			The following example shows that these two notions are independent of each other.
			
			\begin{example}
				Let $(\mathbb{R},d)$ be the upper quasi-metric space (see, Example \ref{e4}). Then $x_n=(n+1)$ is forward parallel to $y_n=(n)$ but $(x_n)$ is not backward parallel to $(y_n)$. \qed
			\end{example}
			Analogous to uniformly continuous functions between metric spaces, forward uniformly continuous functions between quasi-metric spaces can be characterized in terms of forward parallel sequences as follows: a function $f:(X,d) \rightarrow (Y,\rho)$ is forward uniformly continuous if and only if it preserves forward parallel sequences.
		\section{Left K-Cauchy Regular Functions}
		In this section, we study those functions between quasi-metric spaces that preserve left K-Cauchy sequences. We call such functions left K-Cauchy regular.
		\begin{definition}\normalfont A map $f:(X,d)\rightarrow(Y,\rho)$ between two quasi-metric spaces $(X,d)$ and $(Y,\rho)$ is called \textit{left K-Cauchy regular}, if for any left K-Cauchy sequence $(x_n)$ in $(X,d)$, $(f(x_n))$ is left K-Cauchy in $(Y,\rho)$. 
		\end{definition}
		
		So every forward uniformly continuous function is left K-Cauchy regular. Therefore, like metric spaces the class of left-K-Cauchy regular functions contains the class of forward uniformly continuous functions between quasi-metric spaces. 
		
		If a map $f:(X,d)\rightarrow(Y,\rho)$ between two quasi-metric spaces maps right K-Cauchy sequences to  right K-Cauchy sequences, then it is called a \textit{right K-Cauchy regular map}. Note that a sequence in a quasi-metric space $(X,d)$ is $d^s$-Cauchy if and only if it is both left and right K-Cauchy. In \cite{Otafudu}, the author has studied  the class of maps that preserve $d^s$-Cauchy sequences. 
		We call such maps $d^s$-Cauchy regular. We show that every left K-Cauchy regular map is $d^s$-Cauchy regular (Corollary \ref{LeftKCauchyregularimpliesCR}) but the converse may not be true (Example \ref{CRnotLeftKCauchyregular}). First, we give an example of a map which is left K-Cauchy regular but not right K-Cauchy regular.
		\begin{example}\label{e6}
			Let $(X,d)= (\mathbb{N},d)$ and $(Y,\rho)=(\mathbb{R},\rho$), where $d$ is the upper quasi-metric and $\rho$ is the lower quasi-metric (see, Example \ref{e4}). Let $f:(X,d)\rightarrow (Y,\rho)$ be a map defined as $f(x)=x$. Since every left K-Cauchy sequence in $(X,d)$ is eventually constant, $f$ is a left K-Cauchy regular map. But $f$ is not right K-Cauchy regular as the sequence $x_n=n$ is right K-Cauchy in $(X,d)$ and its image is not right K-Cauchy in $(Y,\rho)$. \qed
		\end{example}

		In order to obtain various equivalent characterizations of left K-Cauchy regular maps, we first introduce the concept of semi-equivalent sequences in a quasi-metric space. The concept of semi-equivalent sequences is motivated from the concept of equivalent sequences in a metric space. 
		\begin{definition}\normalfont Let $(x_n)$ and $(y_n)$ be two sequences in a quasi-metric space $(X,d)$. Then
			\begin{enumerate}
				
				\item 	$(x_n)$ is said to be \textit{left semi-equivalent} to $(y_n)$ denoted as $(x_n)\overset{L}{\sim}  (y_n)$ if for any $ \epsilon >0$, there exists $ n_0\in \mathbb{N}$ such that $d(x_m,y_n)<\epsilon, \forall n_0\leq m\leq n$.
				\item $(x_n)$ is said to be \textit{right semi-equivalent} to $(y_n)$ denoted as $(x_n)\overset{R}{\sim}  (y_n)$ if for any $ \epsilon >0$, there exists $ n_0\in \mathbb{N}$ such that $d(y_m,x_n)<\epsilon, \forall n_0\leq m\leq n$
			\end{enumerate}
			If $(x_n)$ is both left as well as right semi-equivalent to $(y_n)$, then we denote it by $(x_n)\underset{R}{\overset{L}{\sim}}(y_n)$.
		\end{definition}

		Recall that two sequences $(x_n)$ and $(y_n)$ in a metric space $(X,d)$ are called \textit{equivalent} if for any $ \epsilon >0$, there exists $ n_0\in \mathbb{N}$ such that $d(x_m,y_n)<\epsilon$ for all $m, n\geq n_0$.
		It is easy to see that in a quasi-metric space $(X,d)$ if two sequence $(x_n)$ and $(y_n)$ are $d^s$-equivalent, then $(x_n)\underset{R}{\overset{L}{\sim}}(y_n)$. However, the converse need not be true. 
		
		The next example shows that if $(x_n)\underset{R}{\overset{L}{\sim}}(y_n)$, then $(x_n)$ may not be $d^s$-equivalent to $(y_n)$; also if $(x_n)$ is left semi-equivalent to $(y_n)$, then $(x_n)$ may not be right semi-equivalent to $(y_n)$.
		\begin{example}
			Consider the upper quasi-metric space on $\mathbb{R}$ as given in Example \ref{e4}. Then the sequence $x_n=-n$ for all $n\in \mathbb{N}$
			is both left and right semi-equivalent to itself. But it is not $d^s$-equivalent to itself. If we consider the sequences $x_n = 2$ and $y_n = 1$ for all \( n \in \mathbb{N} \). Then \((x_n){\overset{L}{\sim}} (y_n) \), but \((x_n){\overset{R}{\nsim}}(y_n)\). \qed		
		\end{example}
		Following important facts about semi-equivalent sequences are easy to prove in a quasi-metric space.
		\begin{itemize}
			\item $(x_n) \overset{L}{\sim} (x_n) \iff (x_n)$ is left K-Cauchy $\iff$ $(x_n) \overset{R}{\sim} (x_n)$.
			\item $(x_n)\buildrel L\over{\sim}  (y_n) \implies (x_n) \overset{\text{f}}{||} (y_n)$.
			\item $(x_n) \underset{R}{\overset{L}{\sim}} (y_n) \implies (x_n) \underset{b}{\overset{f}{||}} (y_n) $.
			\item $(x_n) \underset{R}{\overset{L}{\sim}} (x)\iff (x_n) \underset{b}{\overset{f}{\rightarrow}} (x) $.
			\item $(x_n) \underset{R}{\overset{L}{\sim}} (y_n) \iff (x_1,y_1,x_2,y_2,.....)$ and $(y_1,x_1,y_2,x_2,.....)$ are left K-Cauchy. 
			
		\end{itemize} 
		\begin{remark}
			In the above facts, we cannot replace left K-Cauchy with right K-Cauchy. For instance in the upper quasi-metric space on $\mathbb{R}$, the sequence $x_n=n$ for all $n\in \mathbb{N}$ is right K-Cauchy but $(x_n) \overset{L}{\nsim} (x_n)$ and $(x_n) \overset{R}{\nsim} (x_n)$. And if we take $x_n=y_n=n$ for all $n\in \mathbb{N}$, then $(x_1,y_1,x_2,y_2,.....)$ and $(y_1,x_1,y_2,x_2,.....)$ are right K-Cauchy but $(x_n) \underset{R}{\overset{L}{\nsim}} (y_n)$.
		\end{remark}
		We now provide a sequential characterization of left K-Cauchy regular functions in terms of semi-equivalent sequences. 
		\begin{theorem}\label{LKC characterization}
			Let $(X,d)$ and $(Y,\rho)$ be quasi-metric spaces and $f:(X,d)\rightarrow (Y,\rho)$ be a function. Then the following statements are equivalent:
			\begin{enumerate}[(a)]
				\item $f$ is left K-Cauchy regular;
				\item  if $(x_n)\underset{R}{\overset{L}{\sim}}(y_n)$, then $(f(x_n))\underset{R}{\overset{L}{\sim}}(f(y_n))$.					
			\end{enumerate}
		\end{theorem}
		
		\begin{proof}
			$(a)\implies(b)$. Let $(x_n)\underset{R}{\overset{L}{\sim}}(y_n)$ in $(X,d)$. Then it follows that the sequences $(x_1,y_1,x_2,y_2,\ldots)$ and $(y_1,x_1,y_2,x_2,\ldots)$ are left K-Cauchy sequences in $(X,d)$. Since $f$ is left K-Cauchy regular, $(f(x_1),f(y_1),f(x_2),\ldots)$ and $(f(y_1),f(x_1),f(y_2),\ldots)$ are left K-Cauchy sequences in $(Y,\rho)$. Then it can be seen that for $\epsilon>0,$ there exists $  n_0 \in \mathbb{N}$ such that  $\rho(f(x_m),f(y_n))<\epsilon$   and  $\rho(f(y_m),f(x_n))<\epsilon$ for all $n_0\leq m\leq n$, that is, $(f(x_n))\underset{R}{\overset{L}{\sim}}(f(y_n))$.

			$(b)\implies(a)$. Let $(x_n)$ be a left K-Cauchy sequence in $X$, so $(x_n)$ is both left and right semi-equivalent to itself. Hence $(f(x_n))$ is also both left and right semi-equivalent to itself. Therefore,  $(f(x_n))$ is left K-Cauchy and $f$ is a left K-Cauchy regular map.
		\end{proof}
		
		\begin{proposition}\label{necessary for left K-Cauchy regular}
			Let $f:(X,d)\rightarrow (Y,\rho)$ be a map between quasi-metric spaces $(X,d)$ and $(Y,\rho)$. If $(x_n){\overset{L}{\sim}}(y_n)$ implies $(f(x_n)){\overset{L}{\sim}}(f(y_n))$, then $f$ is left K-Cauchy regular.	
		\end{proposition}
		\begin{proof}
			Suppose $f$ is not left K-Cauchy regular. Therefore, there exists a left K-Cauchy sequence $(x_n)$ in $(X,d)$ such that $(f(x_n))$ is not left K-Cauchy in $(Y,\rho)$. Hence there exists $\epsilon >0$ such that for every $k \in \mathbb{N}$, there exist $n_k, m_k \in \mathbb{N}$ satisfying $n_k\geq m_k > k$ and $\rho (f(x_{m_{k}}),f(x_{n_{k}}))\geq\epsilon$. We can assume that $1 < m_1 \leq n_1 < m_2\leq n_2< \ldots$. Since $(x_n)$ is left K-Cauchy sequence in $(X,d)$, the sequence $(x_{m_{1}},x_{m_{2}},\ldots$) is left semi-equivalent to $(x_{n_{1}},x_{n_{2}},\ldots$). But $(f(x_{m_k}))$ is not left semi-equivalent to $(f(x_{n_k}))$. 
		\end{proof}
		Our next example shows that the converse of the above proposition does not hold. 
		\begin{example}\label{example leftright semiequivalent is must}
			In Example \ref{e6}, the map $f$ is left K-Cauchy regular and the sequence $x_n=2$ for all $n\in \mathbb{N}$ is left semi-equivalent to $y_n=1$ for all $n\in \mathbb{N}$ in $(X,d)$. But $f(x_n)=2$ for all $n\in \mathbb{N}$ is not forward parallel to $f(y_n)=1$ for all $n\in \mathbb{N}$ in $(Y,\rho)$. This example also shows that $f$ is not forward uniformly continuous. \qed
		\end{example}

		\begin{remark}
			Proposition \ref{necessary for left K-Cauchy regular} and Example \ref{example leftright semiequivalent is must} together show that we cannot replace left and right semi-equivalent with left semi-equivalent in Theorem \ref{LKC characterization}.
		\end{remark}
		In order to develop the analysis of left-K-Cauchy regular maps further, we first need to recall the definitions of hereditarily precompact and related concepts. We also need a version of Efremovi$\acute{c}$ lemma in the setting of quasi-metric spaces.

		A subset $Y$ of a quasi-metric space $(X,d)$ is said to be \textit{precompact} (resp. \textit{outside precompact}) if for every $\epsilon>0$ there exists a finite subset $Z$ of $Y$ (resp. of $X$) such that $ Y\subseteq \bigcup \{B^{+}(z,\epsilon):z\in Z\}.$ A subset $Y$ of a quasi-metric space $(X,d)$ is called \textit{hereditarily precompact} if its every subset is precompact. Clearly, 
		\begin{center}
			hereditarily precompact $\Rightarrow$ precompact $\Rightarrow$ outside precompact.
		\end{center} However, reverse implications fail. Also a quasi-metric space is hereditarily precompact if and only if every sequence has a left K-Cauchy subsequence (Proposition 1.2.36, \cite{Cobzas}). For more details about these notions, see \cite{Cobzas,kunzireilly,lambrinos}. Note that in metric spaces all these notions are equivalent and known as total boundedness. Here, we would like to mention that every $d^s$-totally bounded subset of $X$ is hereditarily precompact in $(X,d)$. However, the converse is not true in general. 
		
		Recall that if $d$ is a metric on $X$, then the Efremovi$\acute{c}$ Lemma (page 92, \cite{BEER_BOOK}) states that if $(x_n)$ and $(y_n)$ are sequences in $(X,d)$ and $\epsilon > 0$ such that $d(x_n,y_n)>\epsilon$ for every $n\in \mathbb{N}$, then there exists an infinite subset $E$ of $\mathbb{N}$ such that $d(x_k,y_l)\geq \frac{\epsilon}{4}$ for every $k,l\in E$. 
		
		In a quasi-metric space $(X,d)$, we have $d^s(x,y) \geq d(x,y)$ for all $x,y \in X$. So we have the following version of Efremovi$\acute{c}$ Lemma for quasi-metric spaces. 
		
		\begin{lemma}\label{l2}
			Let $(X,d)$ be a quasi-metric space and $(x_n)$ and $(y_n)$ be sequences in $(X,d)$ such that $d(x_n,y_n)>\epsilon$ for every $n\in \mathbb{N}$. Then there exists an infinite subset $E$ of $\mathbb{N}$ such that $d^{s}(x_k,y_l)\geq \frac{\epsilon}{4}$ for every $k,l\in E$. 
		\end{lemma}
		Our next example shows that in the above lemma, we can not replace $d^s$ by $d$. 
		\begin{example}
			Let $(\mathbb{R},d)$ be the Sorgenfrey quasi-metric space (see, Example \ref{e3}). Consider $x_n= \frac{1}{n}$ and $y_n= \frac{1}{n+1}$ for every $n\in \mathbb{N}$. Now take $\epsilon$ to be any number strictly between $0$ and $1$, then $d(x_n,y_n)>\epsilon$ for every $n\in \mathbb{N}$. Now our claim is that there does not exists any infinite subset $E$ of $\mathbb{N}$ such that $d(x_k,y_l)\geq \frac{\epsilon}{4}$ for every $k,l\in E$. On the contrary, let such an infinite subset $E$ exists. Note that $(x_{n+2})$ is forward parallel to $(y_n)$. Therefore, there exists $n_0\in \mathbb{N}$ such that $d(x_{n+2},y_n)<\frac{\epsilon}{4}$ for every $n\geq n_0$. Since $E$ is an infinite subset of $\mathbb{N}$, there exists $m,l\in E$ such that $l=m+2$ and $d(x_m,y_l)< \frac{\epsilon}{4}$.	\qed
		\end{example}
		
			Recall that in a metric space $(X,d)$, the gap between two nonempty subsets $A$ and $B$ of $X$ is defined as $d(A,B)= \inf \{d(a,b):a\in A \text{ and } b\in B\}$. In \cite{borsik2}, Borsik proved that a map $f: (X,d) \to (Y,\rho)$ between two metric spaces is Cauchy regular if and only if for any two nonempty totally bounded sets $A$ and $B$ in $(X,d)$ such that the gap $d(A,B) = 0$, we have $\rho(f(A),f(B)) = 0$. This is no longer true if $d$ and $\rho$ are assumed to be quasi metrics even if $A$ and $B$ are hereditarily precompact in $(X,d)$. However, we have the following nice result. 
			\begin{theorem}\label{distance_leftKCauchy}
				Let $f: (X,d) \to (Y,\rho)$ be a map between quasi-metric spaces $(X,d)$ and $(Y,\rho)$. Then the following assertions are equivalent:
				\begin{enumerate}[(a)]
					\item $f$ is left K-Cauchy regular; 
					\item if $A$ is hereditarily precompact in $(X,d)$ and $B\subseteq X$ such that $d^{s}(A,B)=0$, then $\rho^{s}(f(A),f(B))=0$;
					\item if $A$ and $B$ are hereditarily precompact sets in $(X,d)$ such that $d^{s}(A,B)=0$, then $\rho^{s}(f(A),f(B))=0$.
				\end{enumerate}
			\end{theorem}
			\begin{proof}
				$(a)\implies (b)$.	Let $A$ be hereditarily precompact subsets of $X$ such that $d^{s}(A,B)=0$. So there exist sequences $(x_n)$ in $A$ and $(y_n)$ in $B$ such that $d^{s}(x_n,y_n)<\frac{1}{n}$. Consequently, $(x_n)\underset{b}{\overset{f}{||}} (y_n)$. Since every sequence in a hereditarily precompact set has a left K-Cauchy subsequence, we can find a left K-Cauchy subsequence $(x_{n_{k}})$ of $(x_n)$ such that $( x_{n_{k}})\underset{R} {\overset{L}{\sim}}( y_{n_{k}})$. Then Theorem \ref{LKC characterization} implies  $ (f(x_{n_{k}})) \underset{R}{\overset{L}{\sim}} (f(y_{n_{k}}))$. Thus, $\rho^{s}(f(A),f(B))=0$.
				$(b)\implies(c)$. It is immediate.
				$(c)\implies(a)$. Suppose $f$ is not left K-Cauchy regular. Then there exist sequences $(x_n)$ and $(y_n)$ in $X$ such that $(x_n)\underset{R} {\overset{L}{\sim}}(y_n)$ but $(f(x_n))\underset{R} {\overset{L}{\nsim}}(f(y_n))$. Without loss of generality, let $(f(x_n)){\overset{L}{\nsim}}(f(y_n))$. Therefore, there exists $\epsilon >0$ such that for every $k \in \mathbb{N}$, there exist $n_k, m_k \in \mathbb{N}$ satisfying $n_k\geq m_k > k$ and $\rho (f(x_{m_{k}}),f(y_{n_{k}}))\geq\epsilon$. We can assume that $1 < m_1 \leq n_1 < m_2\leq n_2< \ldots$. Now from Lemma \ref{l2},  there exists an infinite subset $E$ of $\mathbb{N}$ such that $\rho^{s}(f(x_{m_{k}}),f(y_{n_{l}}))\geq\frac{\epsilon}{4}$ for every $k,l\in E$. Consider $A=\{x_{m_{k}}:k\in E\}$ and $B=\{y_{n_{l}}:l\in E\}$. Since $(x_n)\underset{R} {\overset{L}{\sim}}(y_n)$, the sequences $(x_n)$ and $(y_n)$ are left K-Cauchy sequences in $(X,d)$ and $(x_{m_k})\underset{R} {\overset{L}{\sim}}(y_{n_k})$. Therefore,  $A$ and $B$ are hereditarily precompact subsets of $X$ and $d^{s}(A,B)=0$ but $\rho^{s}(f(A),f(B))\neq 0$.	
			\end{proof}
			In our next example we show that hereditarily precompact cannot be replaced with precompact in Theorem \ref{distance_leftKCauchy}.
			
			\begin{example}\label{e5}
				Let $X=\ell_{\infty}$ (the set of all bounded real sequences) with the quasi-metric $d$ defined as $d(x,y)= \sup \{0,y_1-x_1,y_2-x_2,\dots\}$ for $x, y \in X$. Then the associated metric $d^s$ is the supremum metric on $\ell_{\infty}$, that is, $d^s(x,y)= \sup \{|x_n-y_n|: n \in \mathbb{N}\}$. Now for each $n\in\mathbb{N}$, consider $x_n=\underbrace{(1,1,1,\ldots,1}_\text{n terms},0,0,\dots)$ and $y_n=\underbrace{(1,1,1,\ldots, 1, 1-\frac{1}{n+1}}_\text{n terms},0,0,\dots)$. Then $(x_n)$ forward converges to $x=(1,1,\dots)$ in $(X,d)$. 
				
				Let $A= \{x_n:n\in\mathbb{N}\}\cup \{x\}$ and $B= \{y_n:n\in \mathbb{N}\}$. Then $A$ is precompact but not hereditarily precompact (see, Remark 1.2.24, \cite{Cobzas}).
				
				Define $f:(A\cup B,d)\rightarrow (\mathbb{R},\rho)$, where $\rho$ is the upper quasi-metric, as

				$$f(z)=\Biggl\{\begin{array}{lc}
					2n, \text{   if } z=x_n \text{ for some } n\in \mathbb{N} \\
					2n+1, \text{   if } z=y_n \text{ for some } n\in \mathbb{N} \\
					0, \text{   if } z=x.\end{array}$$ One can easily see that only eventually constant sequences in $A\cup B$ are left K-Cauchy sequence. Hence $f$ is left K-Cauchy regular. But $d^{s}(A,B)=0$ and $\rho^{s}(f(A),f(B))=1$.\qed
			\end{example}

			Recall that for a quasi-metric space $(X,d)$ and $A\subseteq X$, if $A$ totally bounded in $(X,d^s)$, then $A$ is hereditarily precompact in $(X,d)$. Consequently, by Theorem 4 of \cite{borsik2}, we have the following corollary to Theorem \ref{distance_leftKCauchy}.
			\begin{corollary}\label{LeftKCauchyregularimpliesCR}
				Let $(X,d)$ and $(Y,\rho)$ be quasi metric spaces. If $f:(X,d)\rightarrow (Y,\rho)$ is left K-Cauchy regular, then $f:(X,d^s)\rightarrow (Y,\rho^s)$ is a Cauchy regular map.
			\end{corollary}
			
			Our next example shows that the converse of the above corollary need not hold.
			\begin{example}\label{CRnotLeftKCauchyregular}
				Let $(\mathbb{R},d)$ and $(\mathbb{R},\rho)$ be quasi-metric spaces where $d$ is the upper quasi-metric and $\rho$ is the lower quasi-metric. Then $d^s$ and $\rho^s$ be the usual metric on $\mathbb{R}$. Then the identity map $Id:(\mathbb{R},d^s)\rightarrow (\mathbb{R},\rho^{s})$ is Cauchy regular. However, $Id:(\mathbb{R},d)\rightarrow (\mathbb{R},\rho)$ is not left K-Cauchy regular as $x_n=-n$ is left K-Cauchy in $(\mathbb{R},d)$, but $f(x_n)=-n $ is not left K-Cauchy in $(\mathbb{R},\rho)$. \qed
			\end{example}
			
			Beer and Garrido (Theorem 3.2, \cite{beergarrido2014}) has proved that a metric space is totally bounded if and only if each real-valued Cauchy regular map defined on it is bounded. The next theorem is a successful attempt to generalize this result for quasi-metric spaces. First we define bounded subset of a quasi-metric space.
			
			A subset $B$ of a quasi-metric space $(X,d)$ is called $d$-bounded if there exists $x\in X$ and $r>0$ such that $B\subseteq B^+(x,r)$. A function $f:(X,d)\rightarrow (Y,\rho)$ between two quasi-metric spaces is said to be bounded if $f(X)$ is $\rho$-bounded in $(Y,\rho)$.
			
			\begin{theorem}\label{lKC_d-bouded}
				Let $B$ be a nonempty subset of a quasi-metric space $(X,d)$. Then the following are equivalent:
				\begin{enumerate}[(a)]
					\item $B$ is hereditarily precompact;
					\item whenever $(Y,\rho)$ is a quasi-metric space and $f:(X,d)\rightarrow (Y,\rho)$ is left K-Cauchy regular, then $f(B)$ is hereditarily precompact in $(Y,\rho)$;
					\item whenever $(Y,\rho)$ is a quasi-metric space and $f:(X,d)\rightarrow (Y,\rho)$ is left K-Cauchy regular, then $f(B)$ is $\rho$-bounded in $(Y,\rho)$;
					\item whenever $f:(X,d)\rightarrow (\mathbb{R},\rho)$ is left K-Cauchy regular where $\rho$ is an upper quasi-metric on $\mathbb{R}$, then $f(B)$ is bounded in $(\mathbb{R},\rho)$.
				\end{enumerate}
			\end{theorem}
			\begin{proof}
				$(a)\implies (b)$.  Suppose $(y_n)$ be any sequence in $f(B)$. Choose $(x_n)$ in $B$ such that $f(x_n)=y_n$. Since $B$ is hereditarily precompact and $f$ is left K-Cauchy regular, there exists a left K-Cauchy subsequence $(x_{n_{k}})$ in $(x_n)$ and $f(x_{n_{k}})$ is a left K-Cauchy subsequence of $(y_n)$. Hence $f(B)$ is hereditarily precompact in $(Y,\rho)$.
				
				$(b)\implies (c)$. It follows from the fact that every hereditarily precompact subset of a quasi-metric space $(X,d)$ is $d$-bounded.

				$(c)\implies (d)$. It is immediate.

				$(d)\implies(a)$. On contrary, suppose $B$ is not hereditarily precompact.  Then there exists an infinite subset $A \subseteq B$ which is not precompact. Hence there exists $\epsilon > 0$ such that 
				$$A \nsubseteq \bigcup_{x \in F} B^{+}(x,\epsilon)$$
				for any finite subset $F \subseteq A$.  
				
				Choose an element $x_1 \in A$. Since $A \nsubseteq B^{+}(x_1,\epsilon)$, we can choose $x_2 \in A$ such that $x_2 \notin B^{+}(x_1,\epsilon)$. Similarly, choose $x_3 \notin B^{+}(x_1,\epsilon)\cup B^{+}(x_2,\epsilon)$, and so on. In this way we obtain a sequence $(x_n)$ in $A$ such that 
				$$d(x_k,x_{k+i})>\epsilon \quad \text{for all } i,k\in\mathbb{N}.$$
				
				Now set 
				$$\delta_{1}=\min\{\epsilon,1\}, \quad 
				\delta_{2}=\min\{\epsilon,\tfrac{1}{2},d(x_2,x_1)\}, \quad 
				\delta_{3}=\min\{\epsilon,\tfrac{1}{3},d(x_3,x_1),d(x_3,x_2),\delta_{2}\},$$
				and continue similarly to define $\delta_n$ for each $n \in \mathbb{N}$.
				
				Define a map $f:(X,d)\to (\mathbb{R},\rho)$, where $\rho$ is the upper quasi-metric defined in Example~\ref{e4}, by
				$$
				f(x)=
				\begin{cases}
					n-\dfrac{4n}{\delta_{n}}\, d(x,x_n), & \text{if } x \in B^{d^{s}}(x_n,\tfrac{\delta_{n}}{4}) \text{ where $n$ is the least such natural number}, \\[6pt]
					0, & \text{otherwise}.
				\end{cases}
				$$
				
				Before showing that $f$ is left $K$-Cauchy regular, note that if $(y_n)$ is a left $K$-Cauchy sequence in $(X,d)$, then $(y_n)$ cannot contain more than one element of $(x_n)$ occurring infinitely many times, since $d(x_k,x_{k+i})>\epsilon$ for all $i,k\in \mathbb{N}$.
				
				We now consider the following cases.

				\noindent\textbf{Case 1.} Suppose $(y_n)$ is a left $K$-Cauchy sequence and some $x_k \in (x_n)$ occurs infinitely often in $(y_n)$. Since $(y_n)$ is left $K$-Cauchy, 
				$$y_n \underset{b}{\overset{f}{\rightarrow}} x_k,$$
				which implies that there exists $n_0 \in \mathbb{N}$ such that $y_n \in B^{d^{s}}(x_k,\tfrac{\delta_{k}}{4})$ for every $n \geq n_0$. Hence $f(y_n)\xrightarrow{d^{s}} x_k$, so $(f(y_n))$ is $\rho^{s}$-Cauchy, and therefore left $K$-Cauchy.  
				
				Moreover, if infinitely many $y_n \in B^{d^{s}}(x_i,\tfrac{\delta_{i}}{4})$ for some $i<k$, then $d(x_i,x_k)<\epsilon$, contradicting the construction of $(x_n)$.

				\noindent\textbf{Case 2.} Suppose $(y_n)$ is left $K$-Cauchy such that for infinitely many $n$ (say $n_k$) we have $y_{n_k}\in B^{d^{s}}(x_{i_k},\tfrac{\delta_{n}}{4})$. Then
				$$(y_{n_k}) \underset{b}{\overset{f}{||}} (x_{i_k}),$$
				so that 
				$$d(x_{i_k},x_{i_l}) \leq d(x_{i_k},y_{n_k}) + d(y_{n_k},y_{n_l}) < \epsilon \quad \text{for } l > K \text{ sufficiently large},$$
				which is impossible.

				\noindent\textbf{Case 3.} Suppose $(y_n)$ is left $K$-Cauchy and infinitely many terms of $(y_n)$ belong to finitely many such balls.

				\emph{Subcase 3a.} If infinitely many terms of $(y_n)$ belong to both $B^{d^{s}}(x_k,\tfrac{\delta_{k}}{4})$ and $B^{d^{s}}(x_l,\tfrac{\delta_{l}}{4})$, then by left $K$-Cauchiness we obtain $d(x_k,x_l)<\epsilon$, a contradiction.

				\emph{Subcase 3b.} If infinitely many terms of $(y_n)$ belong to only one ball, say $B^{d^{s}}(x_l,\tfrac{\delta_{l}}{4})$, and infinitely many other terms lie outside all such balls, then we obtain a subsequence $(y_{m_k})$ with $d^{s}(x_l,y_{m_k})>\tfrac{\delta_{l}}{4}$. For every $\gamma$ with $0<\gamma<\tfrac{\delta_{l}}{4}$, we get $d^{s}(x_l,y_{n_{k}})>\gamma$ eventually. Hence $d^{s}(x_l,y_{n_k})\to \tfrac{\delta_{l}}{4}$, which implies $f(y_{n_k})\xrightarrow{\rho^{s}} 0$ while $f(y_{m_k})=0$ for all $k\in\mathbb{N}$. Thus $(f(y_n))$ is $\rho^s$ Cauchy and hence left $K$-Cauchy in $(\mathbb{R},\rho)$.

				\emph{Subcase 3c.} If all but finitely many terms of $(y_n)$ belong to a single ball, say $B^{d^{s}}(x_k,\tfrac{\delta_{k}}{4})$, then
				$$f(y_n)-f(y_m)=k-\frac{4}{\delta_{k}}d(y_n,x_k)-k+\frac{4}{\delta_{k}}d(y_m,x_k)\leq \frac{4}{\delta_{k}}d(y_m,y_n),$$
				which shows that $(f(y_n))$ is left $K$-Cauchy in $(\mathbb{R},\rho)$.

				In all cases, we conclude that $f$ is left $K$-Cauchy regular.
			\end{proof}

			The above theorem shows that a left K-Cauchy regular function maps hereditarily precompact sets to hereditarily precompact sets. In the context of metric spaces, such maps were first considered by Beer and Levi in \cite{Beer_totalBounded}. Later, Lipsy and Kundu studied these maps extensively in \cite{Lipsy}, where they called them Cauchy subregular maps.

			Our next example illustrates that Theorem \ref{lKC_d-bouded} may fail if we replace hereditarily precompact with precompact.
			\begin{example}
				In Example \ref{e5}, $f$ is left K-Cauchy regular and $A$ is precompact. But $f(A)$ is not bounded in $(\mathbb{R},\rho)$. \qed
			\end{example}

			Note that Theorem \ref{lKC_d-bouded} still holds if we consider lower quasi-metric on $\mathbb{R}$ in place of upper quasi-metric. The next corollary shows that under the assumption of left K-completeness on $(X,d)$, the hereditarily precompatness of a subset of $(X,d)$ can be characterized in terms of boundedness of its images under the left K-Cauchy regular maps defined from $(X,d)$ to $(\mathbb{R}, | \cdot |)$.

			\begin{corollary}
				Let $B$ be a nonempty subset of a left K-complete quasi-metric space $(X,d)$. Then the following are equivalent:
				\begin{enumerate}[(a)]
					\item $B$ is hereditarily precompact;
					\item whenever $(Y,\rho)$ is a quasi-metric space and $f:(X,d)\rightarrow (Y,\rho)$ is left K-Cauchy regular, then $f(B)$ is hereditarily precompact in $(Y,\rho)$;
					\item whenever $(Y,\rho)$ is a quasi-metric space and $f:(X,d)\rightarrow (Y,\rho)$ is left K-Cauchy regular, then $f(B)$ is $\rho$-bounded in $(Y,\rho)$;
					\item whenever $f:(X,d)\rightarrow (\mathbb{R},|\cdot|)$ is left K-Cauchy regular, then $f(B)$ is bounded in $(\mathbb{R},|\cdot|)$.
				\end{enumerate}
			\end{corollary}
			\begin{proof}
				$(a)\implies (b)\implies (c)\implies (d)$ follows from Theorem \ref{lKC_d-bouded}.\\
				$(d)\implies (a)$  On contrary, suppose $B$ is not hereditarily precompact.  Then there exists an infinite subset $A \subseteq B$ which is not precompact. Construct the sequence $(x_n)$ in $A$ and $\delta_n$ for each $n\in\mathbb{N}$ in the manner similar to $(d)\implies (a)$ in Theorem \ref{lKC_d-bouded}.\\
				Now define a map $f:(X,d)\to (\mathbb{R},|\cdot|)$ by
				$$
				f(x)=
				\begin{cases}
					n-\dfrac{4n}{\delta_{n}}\, d(x_n,x), & \text{if } x \in B^{d^{s}}(x_n,\tfrac{\delta_{n}}{4}) \text{ where $n$ is the least such natural number}, \\[6pt]
					0, & \text{otherwise}.
				\end{cases}$$
				Note that Cases $1,2,3a,3b$ can be proved in the similar manner as in $(d)\implies (a)$ of Theorem~\ref{lKC_d-bouded}, except for the Case $3c$, that is, suppose all but finitely many terms of $(y_n)$ belong to a single ball, say $B^{d^{s}}(x_k,\tfrac{\delta_{k}}{4})$. Since $(y_n)$ is left $K$-Cauchy in $(X,d)$ and $(X,d)$ is left $K$-complete, there exists $a\in X$ such that $y_n\xrightarrow{f} a$. Let $d(x_k,a)=p$. Then, as $d(x_k,y_n)\leq d(x_k,a)+d(a,y_n)$, we obtain $d(x_k,y_n)\to p$. Consequently, $f(y_n)\to k-\tfrac{4p}{\delta_{k}}$, which shows that $(f(y_n))$ converges in $(\mathbb{R},|\cdot|)$ and hence $(f(y_n))$ is Cauchy in $(\mathbb{R},|\cdot|)$.
			\end{proof}

			Our last result of this section says that the collection of all left K-Cauchy regular maps from a quasi-metric space $(X,d)$ to $(\mathbb{R},\rho)$ forms a lattice cone, where $\rho$ is the upper quasi-metric on $\mathbb{R}$. A lattice cone $C$ is a cone that carries a partial order relation $\leq$ such that (a) $x\leq y$ implies $x+z\leq y+z$ and $\alpha x\leq \alpha y$ for all $x,y,z\in C$ and for every nonnegative real number $\alpha$; and (b) supremum and infimum of every pair of elements is again an element of the cone \cite{WalterRoth}.
			\begin{theorem}\label{algebraic_properties_of_LKC}
				Let $(X,d)$ be a quasi-metric space and $\rho$ denotes the upper quasi-metric on $\mathbb{R}$. Then the collection
				$$LKC(X)=\{f\in \mathbb{R}^X: f:(X,d)\rightarrow (\mathbb{R},\rho) \text{ is left K-Cauchy regular}\}$$
				
				\begin{enumerate}[(a)]
					\item forms a lattice cone.
					
					\item is closed under the pointwise product of non-negative elements.
				\end{enumerate}
			\end{theorem}
			
			\begin{proof}
				$(a)$. It is easy to see for $c\geq 0$ and $f\in LKC(X)$, we have $cf\in LKC(X)$. 
				
				Now let $f,g\in LKC(X)$, $(x_n)$ be a left K-Cauchy sequence in $(X,d)$ and $\epsilon >0$. Then there exists $n_0 \in \mathbb{N}$ such that $\rho(f(x_m), f(x_n)) < \frac{\epsilon}{2}$ and $\rho(g(x_m), g(x_n)) < \frac{\epsilon}{2}$ for $n\geq m 
				\geq n_0$. As\\

				$\rho(f(x_m)+g(x_m),f(x_n)+g(x_n)) = $
				\[\begin{cases} 
					
					f(x_n)+g(x_n) - f(x_m)-g(x_m), & \text{if } f(x_n)+g(x_n) \geq f(x_m)+g(x_m)\\
					0, & \text{if } f(x_n)+g(x_n) < f(x_m)+g(x_m).
				\end{cases}
				\]

				Therefore, $\rho(f(x_m)+g(x_m),f(x_n)+g(x_n)) \leq \rho (f(x_m),f(x_n))+\rho(g(x_m),g(x_n))$. Hence for $n\geq m\geq n_0$, we get $\rho(f(x_m)+g(x_m),f(x_n)+g(x_n)) < \epsilon$. Hence $f+g\in LKC(X)$.
				
				Now we show that for $f,g \in LKC(X)$, the map $h(x)= \max\{f(x),g(x)\}\in LKC(X)$.
				
				Let $(x_n)$ be a left $K$-Cauchy sequence in $(X,d)$ and fix $\epsilon > 0$. Since $(f(x_n))$ and $(g(x_n))$ are left $K$-Cauchy sequences, there exists $n_0 \in \mathbb{N}$ such that $\rho(f(x_n),f(x_m)) < \epsilon$ and $\rho(g(x_n),g(x_m)) < \epsilon$ for all $m \geq n \geq n_0$. We now analyze the non-trivial cases of $h(x_n)$ and $h(x_m)$.

				\noindent
				\textbf{Case 1.} Suppose $h(x_n) = f(x_n)$ and $h(x_m) = g(x_m)$. Then  
				$\rho(h(x_n),h(x_m)) = \rho(f(x_n),g(x_m))$ $= g(x_m) - f(x_n)$ if $g(x_m) \geq f(x_n)$, and $\rho(f(x_n),g(x_m)) = 0$ otherwise. Since $g(x_m) \geq f(x_n) \geq g(x_n)$, we obtain $g(x_m) - f(x_n) \leq g(x_m) - g(x_n)\leq \rho(g(x_n),g(x_m)) $.

				\noindent
				\textbf{Case 2.} Suppose $h(x_n) = g(x_n)$ and $h(x_m) = f(x_m)$. Then  
				$\rho(h(x_n),h(x_m))= \rho(g(x_n),f(x_m))$ $= f(x_m) - g(x_n)$ if $f(x_m) \geq g(x_n)$, and $\rho(g(x_n),f(x_m)) = 0$ otherwise. Since $f(x_m) \geq g(x_n) \geq f(x_n)$, we obtain $f(x_m) - g(x_n) \leq f(x_m) - f(x_n)\leq \rho(f(x_n),f(x_m))$.

				In either case, the required inequality holds, which shows that $h \in LKC(X)$. Finally, by an entirely analogous argument, one verifies that the pointwise minimum of two left $K$-Cauchy regular functions is again a left $K$-Cauchy regular map. 
				
				$(b)$.	Let $(x_n)$ be a left K-Cauchy sequence in $(X,d)$ and $\epsilon >0$. Since $f,g \in LKC(X)$, $(f(x_n))$ and $(g(x_n))$ are left K-Cauchy sequences in $(\mathbb{R},\rho)$. So we can find $k> 0$ and $n_0 \in \mathbb{N}$ such that $f(x_n)\leq k$ and $g(x_n)\leq k$ for all $n\in \mathbb{N}$,  and $\rho(f(x_m), f(x_n)) < \frac{\epsilon}{2k}$ and $\rho(g(x_m), g(x_n)) < \frac{\epsilon}{2k}$ for $n\geq m 
				\geq n_0$. Now\\

				$\rho(f(x_m)g(x_m),f(x_n)g(x_n)) =$\[\begin{cases}
					0, & \text{   if } f(x_n)g(x_n) < f(x_m)g(x_m)\\
					f(x_n)g(x_n) - f(x_m)g(x_m), & \text{   if }f(x_n)g(x_n) \geq f(x_m)g(x_m). 
				\end{cases}\]
				Since $f,g\geq 0$, it follows that $$\begin{aligned}
					f(x_n)g(x_n) - f(x_m)g(x_m)&= f(x_n)g(x_n)-f(x_n)g(x_m)+f(x_n)g(x_m) - f(x_m)g(x_m)\\ &\leq f(x_n)\rho(g(x_m),g(x_n))+ g(x_m)\rho (f(x_m),f(x_n)).
				\end{aligned}$$ Consequently, $\rho(f(x_m)g(x_m),f(x_n)g(x_n)) < \epsilon$ for all $n\geq m\geq n_0$. Hence $fg \in LKC(X)$.
			\end{proof}

			\begin{remark}
				In Theorem \ref{algebraic_properties_of_LKC},  if we take $\rho$ to be the usual metric on $\mathbb{R}$, then $LKC(X)$ forms a vector lattice and is closed under pointwise product.
				\end{remark}
			
			The following example demonstrates that we cannot drop non-negativity condition in part $(b)$ of Theorem \ref{algebraic_properties_of_LKC}. 
			\begin{example}
				Let $(X,d)=(Y,\rho)=(\mathbb{R},\rho$), where $\rho$ is the upper quasi-metric. Then the maps $f,g:(X,d)\rightarrow (Y,\rho)$ defined as $f(x)=x$ and $g(x)=-1$ are left K-Cauchy regular maps. However, $fg$ is not left K-Cauchy regular as the sequence $x_n=-n$ is left K-Cauchy in $(X,d)$ but $f(x_n)g(x_n)=n$ is not a left K-Cauchy sequence in $(Y,\rho)$.\qed
			\end{example}

			\section{Left K-Cauchy Regular and ff-Continuous Functions}
			It is well-known that every Cauchy regular function on a metric space is continuous. Consequently, these functions are also called Cauchy continuous functions. However, in quasi-metric spaces this correspondence fails. This section is devoted to study the relationship between the left K-Cauchy regular functions and ff-continuous functions.  We first give an example of a function which is left K-Cauchy regular but not ff-continuous.
			\begin{example}\label{e10}
				Let $X=Y=(\mathbb{R}^+ \cup\{0\})$. Let $d$ be the Sorgenfrey quasi-metric on $X$ and $\rho$ be the usual metric on $Y$. Then the function $f:(X,d)\rightarrow (Y,\rho)$ defined as $$f(x)=\Biggl\{\begin{array}{lc}
					\frac{1}{x}, \text{   if } x\neq 0 \\
					0, \text{   if } x=0
				\end{array}$$ is left K-Cauchy regular as every left K-Cauchy sequence in $(X,d)$ is eventually increasing. However, $f$ is not ff-continuous at $x=0$ as the sequence $(1/n)$ is forward convergent to $0$ in $(X,d)$ but $(f(\frac{1}{n}))$ is not forward convergent to $f(0)$ in $(Y,\rho)$.	\qed
			\end{example}
			
			\begin{theorem}\label{leftKCauchy_implies_ffcts} Let $(X,d)$ and $(Y,\rho)$ be quasi-metric spaces, where $(X,d)$ is $T_1$ and has the property forward convergence implies backward convergence. If $f:(X,d) \rightarrow (Y,\rho)$ is a left K-Cauchy regular map, then $f$ is ff-continuous.
			\end{theorem}
			\begin{proof}
				Let $x\in X$ and $(x_n)$ be a sequence in $X$ such that $x_{n}\xrightarrow{f}x$. Since $(X,d)$ is $T_1$ and has property forward convergence implies backward convergence, the constant sequence $(x)$ is left and right semi-equivalent to $(x_n)$. So by Theorem \ref{LKC characterization}, $(f(x))$ is left and right semi-equivalent to $(f(x_n))$. Hence $f(x_n)\xrightarrow{f}f(x)$. Consequently, $f$ is ff-continuous at $x$. 
			\end{proof}
			Example \ref{e10} shows that above theorem is not true if we drop the condition  forward convergence implies backward convergence as the sequence $(1/n)$ is forward convergent to $0$, but it is not backward convergent in $(X,d)$.

			The next theorem provides a condition under which the converse of the Theorem \ref{leftKCauchy_implies_ffcts} holds. 
			
			\begin{theorem}\label{ffcts_implies_leftKCauchyregular}
				Let $(X,d)$ be a left K-complete quasi-metric space and $(Y,\rho)$ be a $T_1$ quasi-metric space having the property forward convergence implies backward convergence. Then every ff-continuous function from $(X,d)$ to $(Y,\rho)$ is left K-Cauchy regular. 
			\end{theorem}
			
			\begin{proof}
				Let $(x_n)\underset{R}{\overset{L}{\sim}}(y_n)$ in $(X,d)$. Then the sequence $(z_n)=(x_1,y_1,x_2,y_2,\ldots)$ is left K-Cauchy in $(X,d)$. Since $(X,d)$ is left K-complete, there exists $x\in X$ such that $z_{n}\xrightarrow{f}x$. So $f(z_n)\xrightarrow{f}f(x)$. As $(Y,\rho)$ is $T_1$ and has property forward convergence implies backward convergence, we can find $n_0 \in \mathbb{N}$ such that $\rho(f(x_m),f(y_n))\leq \rho(f(x_m),f(x))+\rho(f(x),f(y_n))$ and $\rho(f(y_m),f(x_n))\leq \rho(f(y_m),f(x))+\rho(f(x),f(x_n))$ for all $n\geq m\geq n_0$. Hence $(f(x_n))\underset{R}{\overset{L}{\sim}}(f(y_n))$.
			\end{proof}  
			The next example shows that the condition forward convergence implies backward convergence on the space \((Y, \rho)\) cannot be dropped.

			\begin{example}
				Let $X=\{\frac{1}{n}:n\in\mathbb{N}\}\cup \{0\}$ with the usual metric on it and let $(Y,\rho)$ be the space of real numbers with the Sorgenfrey quasi-metric. Note that $(X,|\cdot|)$ is left K-complete and $(Y,\rho)$ does not have the property forward convergence implies backward convergence as the sequence $(1/n)$ is forward convergent to $0$ but not backward convergent in $(Y,\rho)$. Now the map $f:(X,|\cdot|)\rightarrow (Y,\rho)$ defined as $f(x)=x$ is ff-continuous. However, $f$ is not left K-Cauchy regular as $(\frac{1}{n})$ is left K-Cauchy in $(X,|\cdot|)$ but it is not left K-Cauchy in $(Y,\rho)$.\qed
			\end{example}
			
			\begin{corollary}
				Let $(X,d)$ be a left K-complete quasi-metric space. Then every real-valued ff-continuous function from $(X,d)$ to $(\mathbb{R},|\cdot|)$ is left K-Cauchy regular. 
			\end{corollary}
			
			\begin{theorem}
				Let $(X,d)$ and $(Y,\rho)$ be $T_1$ quasi-metric spaces with the property forward convergence implies backward convergence. If $f:(X,d)\rightarrow (Y,\rho)$ be a bijection such that $f$ and $f^{-1}$ are left K-Cauchy regular maps, then $(X,d)$ is left K-complete if and only if $(Y,\rho)$ is left K-complete.
			\end{theorem}
			\begin{proof}
				We only need to prove that if $(X,d)$ is left K-complete, then $(Y,\rho)$ is left K-complete. Let $(y_n)$ be a left K-Cauchy sequence in $(Y,\rho)$. Then $(f^{-1}(y_n))$ is left K-Cauchy in $(X,d)$. Since $(X,d)$ is left K-complete, $(f^{-1}(y_n))$ is forward convergent in $(X,d)$. By Theorem \ref{leftKCauchy_implies_ffcts}, $f$ is ff-continuous. Therefore $(f(f^{-1}(y_n)))=(y_n)$ is forward convergent in $(Y,\rho)$.
			\end{proof}
			
			\begin{theorem}
				Let $(Y,\rho)$ be a quasi metric space. If every left K-Cauchy regular map from any quasi-metric space $(X,d)$ to $(Y,\rho)$ is ff-continuous,
				then $(Y,\rho)$ is left K-Complete.
			\end{theorem}
			\begin{proof}
				Suppose there exists a left K-Cauchy sequence $(y_n)$ in $(Y,\rho)$ which does not converge. Let $X=\{n:n\in \mathbb{N}\}$ be a subspace of $\mathbb{R}$ with lower quasi-metric $d$. Then the map $f:(X,d)\rightarrow(Y,\rho)$ defined as $f(n)=y_n$ is left K-Cauchy regular but not ff-continuous at any $n\in\mathbb{N}$.
			\end{proof}
			
			In our last result of this section, we demonstrate the conditions under which a left K-Cauchy regular map defined on a subset of a quasi-metric space can be extended to its forward closure. The closure of a subset $A$ of a quasi-metric space $(X,d)$ is called the forward closure of $A$ and is denoted by $Cl^+A$.

			\begin{theorem}\label{t3}
				Let $(X,d)$ and $(Y,\rho)$ be $T_1$ quasi-metric spaces with the property forward convergence implies backward convergence. If $(Y,\rho)$ is left K-Complete, $A\subseteq X$ and $f:(A,d)\rightarrow (Y,\rho)$ is a left K-Cauchy regular function, then there is a unique extension $\bar{f}:(Cl^+{A},d)\rightarrow (Y,\rho)$ of $f$ such that $\bar{f}$ is left K-Cauchy regular. 
			\end{theorem}
			\begin{proof}
				Let $x\in Cl^+{A}$. Then there exists a sequence $(x_n)$ in $A$ such that $x_n\xrightarrow{f} x$. Since $(X,d)$ is $T_1$ with the property forward convergence implies backward convergence and $f:(A,d)\rightarrow (Y,\rho)$ is left K-Cauchy regular, $(f(x_n))$ is left K-Cauchy in $(Y,\rho)$. As $(Y,\rho)$ is left K-Complete and $T_1$ with the property that forward convergence implies backward convergence, there exists unique $y\in Y$ such that $f(x_n)\underset{b}{\overset{f}{\rightarrow}} y$. 
				Now suppose $(y_n)$ is another sequence in $(A,d)$ such that $y_n\xrightarrow{f} x$.  Then $(x_n)\underset{R}{{\overset{L}{\sim}}}(y_n)$ and by Theorem \ref{LKC characterization}, $(f(x_n))\underset{R}{{\overset{L}{\sim}}}(f(y_n))$. Hence $f(y_n)\underset{b}{\overset{f}{\rightarrow}} y$. 
				
				Define $\bar{f}:(Cl^+{A},d)\rightarrow (Y,\rho)$ by
				$\bar{f}(x)=y.$ Clearly $\bar{f}(x)=f(x)$ for $x\in A$. Note that by the construction of $\bar{f}$, for any $x \in Cl^+ A$, and $\delta>0, \epsilon>0$, there exists $a\in A$ such that $d(x,a)<\delta$, $d(a,x)<\delta$ and $\rho(\bar{f}(x),f(a))<\epsilon$, $\rho(f(a),\bar{f}(x))<\epsilon$. 
				
				Now we show that $\bar{f}$ is left K-Cauchy regular. Let $(x_n)$ be a left K-Cauchy sequence in $(Cl^+ A,d)$. So for each $n\in \mathbb{N}$, there exists $a_n \in A$ such that $d(x_n,a_n)<1/n$, $d(a_n,x_n)<1/n$ and  $\rho(\bar{f}(x_n),f(a_n))<1/n$, $\rho(f(a_n),\bar{f}(x_n))<1/n$. As $(x_n)$ is a left K-Cauchy sequence in $(Cl^+ A,d)$, $(x_n)\underset{b}{\overset{f}{||}} (a_n)$ and $d(a_n,a_m)\leq d(a_n,x_n)+d(x_n,x_m)+d(x_m,a_m)$, we have $(a_n)$ left K-Cauchy in $(A,d)$. Therefore, $f(a_n)$ is left K-Cauchy in $(Y,\rho)$. Now let $\epsilon >0$. Then there exists $n_0$ such that for all $n\geq m\geq n_0$, $\rho(f(a_m),f(a_n))<\epsilon/3$, $\rho(\bar{f}(x_m),f(a_m))<\epsilon/3$ and $\rho(f(a_n),\bar{f}(x_n))<\epsilon/3$. Hence $\rho(\bar{f}(x_m),\bar{f}(x_n)) \leq \rho(\bar{f}(x_m),f(a_m))+\rho(f(a_m),f(a_n))+\rho(f(a_n),\bar{f}(x_n))<\epsilon/3+ \epsilon/3+ \epsilon/3=\epsilon$ for all $n\geq m\geq n_0$. So $\bar{f}$ is left K-Cauchy regular.
				
				Uniqueness of $\bar{f}$ follows from the Remark \ref{R1} and Theorem \ref{leftKCauchy_implies_ffcts}.
			\end{proof}

			\section{Left K-Cauchy Regular and Forward Uniformly Continuous Functions}
			This section is devoted to study the relationship between the forward uniformly continuous functions and left K-Cauchy regular functions. We have already seen that every forward uniformly continuous function between quasi-metric spaces is always left K-Cauchy regular. However, the converse need not be true even for the maps between metric spaces. In the realm of metric spaces, Snipes (see \cite{snipes3}, Theorem 3), showed that a Cauchy regular map on a totally bounded metric spaces is uniformly continuous. In general, this result fails for left K-Cauchy regular maps between quasi-metric spaces. 
			
			\begin{example}\label{CRBNUC}
				Let $X=\{\frac{1}{n}:n\in\mathbb{N}\}\cup \{0\}$ and $d$ be the conjugate of the Sorgenfrey quasi-metric. Since an eventually decreasing  sequence in $(X,d)$ is left K-Cauchy, $(X,d)$ is hereditarily precompact. Define $f:(X,d)\rightarrow (\mathbb{R},\rho)$, where $\rho$ is the lower quasi-metric, as $$f(x)=\Biggl\{\begin{array}{lc}
					\frac{1}{x}, \text{   if } x\neq 0 \\
					0, \text{   if } x=0.
				\end{array}$$ 
				Then it is easy to see that $f$ is left K-Cauchy regular. However, $f$ is not forward uniformly continuous as $(\frac{1}{n})$ is forward parallel to the constant sequence $(0)$ in $X$ but $(f(x_n))$ is not forward parallel to $(f(y_n))$.\qed
			\end{example}
			
			More generally, a map between metric spaces is Cauchy regular if and only if its restriction to each totally bounded subset is uniformly continuous (see,  \cite{Beer_totalBounded}, \cite{borsik2}). Now, we discuss the conditions under which this result holds in the setting of quasi-metric spaces.
			
			\begin{theorem}\label{t8}
				Let $(X,d)$ and $(Y,\rho)$ be quasi-metric spaces, where $(X,d)$ has the property forward parallel implies backward parallel. Let $f:(X,d)\rightarrow (Y,\rho)$ be a map. Then the following are equivalent:
				\begin{enumerate}[(a)]
					\item $f$ is left K-Cauchy regular;
					\item restriction of $f$ to each non-empty hereditarily precompact subset of $X$ is forward uniformly continuous;
					\item restriction of $f$ to each non-empty hereditarily precompact subset of $X$ is left K-Cauchy regular.
				\end{enumerate}
			\end{theorem}
			\begin{proof}
				$(a)\implies (b)$.	Suppose that there exists a hereditarily precompact subset $A$ of $X$ such that $f_{|A}:X\rightarrow Y$ is not forward uniformly continuous. Then there exists sequences $(x_n)$ and $(y_n)$ in $A$ such that $ (x_n)\overset{f}{||}(y_n)$ in $X$ but $(f(x_n))\overset{f}{\not||}(f(y_n))$ in $Y$. Therefore, there exists $\epsilon >0$ such that for every $k\in \mathbb{N}$, there exists $n_k \in\mathbb{N}$ such that $\rho (f(x_{n_{k}}),f(y_{n_{k}}))>\epsilon$. As $A$ is hereditarily precompact, there exists left-K-Cauchy subsequences $(x_{n_{k_{j_{m}}}})$ and $(y_{n_{k_{j_{m}}}})$ of $(x_n)$ and $(y_n)$ in $X$. Since $(X,d)$ has the property forward parallel implies backward parallel, $(x_{n_{k_{j_{m}}}})\underset{R}{{\overset{L}{\sim}}}(y_{n_{k_{j_{m}}}})$. By Theorem \ref{LKC characterization}, $f(x_{n_{k_{j_{m}}}})\underset{b}{\overset{{f}}{||}}f(y_{n_{k_{j_{m}}}})$ in $Y$. So we arrive at a contradiction.
				
				$(b)\implies (c)$. It follows from the fact that every forward uniformly continuous function is left K-Cauchy regular.
				
				$(c)\implies (a)$. It follows from the fact that every left K-Cauchy sequence forms a hereditarily precompact set.
			\end{proof}
			
			\begin{corollary}
				Let $(X,d)$ and $(Y,\rho)$ be quasi-metric spaces. If $(X,d)$ is hereditarily precompact and has the property forward parallel implies backward parallel, then every left K-Cauchy regular map is forward uniformly continuous.
			\end{corollary}

			\begin{example}
				Example \ref{CRBNUC} shows that the property forward parallel implies backward parallel can't be dropped in the above result. As $(\frac{1}{n})\overset{f}{||}(\frac{1}{n+1})$ but $(\frac{1}{n})\overset{b}{\not||}(\frac{1}{n+1})$ in $(X,d)$.\qed
			\end{example}
			
			\begin{remark}
				Note that in a quasi-metric space the property `forward parallel implies backward parallel' is stronger than the property `forward convergence implies backward convergence'. In Example \ref{CRBNUC}, the space $(X,d)$ does not have the property forward parallel implies backward parallel, but it has the property forward convergence implies backward convergence. Therefore, Example \ref{CRBNUC} also demonstrates that the property forward parallel implies backward parallel cannot be replaced with forward convergence implies backward convergence in the above result. 
			\end{remark}

			In \cite{beergariddo2020}, Beer and Garrido explored in detail the phenomenon that a function between metric spaces has a given metric property if and only if its composition with an arbitrary real-valued Lipschitz function has same property. For example,  in Corollary 3.5 of \cite{beergariddo2020}, they showed that a map between metric spaces is Cauchy-regular if and only if its composition with an arbitrary real-valued Lipschitz function is Cauchy-regular. In order to prove a similar result for maps between quasi-metric spaces, we first need the following result. 
			\begin{proposition}\label{funiformlycts}
				Let $(X,d)$ be a quasi-metric space with the property forward parallel implies backward parallel and let $A$ be a nonempty subset of $X$. Then the functions $d(\cdot,A):(X,d)\rightarrow (\mathbb{R},|\cdot|)$  and $d(A,\cdot):(X,d)\rightarrow (\mathbb{R},|\cdot|)$ defined as $$d(x,A)=\inf\{ d(x,a):a\in A\}$$
				$$d(A,x)=\inf\{ d(a,x):a\in A\}$$ are forward uniformly continuous.
			\end{proposition}
			\begin{proof}
				Let $(x_n)$ be forward parallel to $(y_n)$  in $X$. Since $(X,d)$ has the property forward parallel implies backward parallel, we have $(x_n)\underset{b}{\overset{f}{||}} (y_n)$. Also  $d(x_n,A)-d(y_n,A)\leq d(x_n,y_n)$ and  $d(y_n,A)-d(x_n,A)\leq d(y_n,x_n)$. Let $\epsilon > 0$. Since $(x_n)\underset{b}{\overset{f}{||}} (y_n)$, there exists $n_0 \in \mathbb{N}$ such that $d(x_n, y_n) < \epsilon$ and $d(y_n, x_n) < \epsilon$. Consequently, we have $|d(x_n,A)-d(y_n,A)| < \epsilon$. Hence $d(\cdot,A)$ is forward uniformly continuous. Similarly, we can show that $d(A,\cdot)$ is forward uniformly continuous.
			\end{proof}
			
			\begin{theorem}
				Let $f:(X,d)\rightarrow (Y,\rho)$ be a map between two quasi-metric spaces $(X,d)$ and $(Y,\rho)$. If $(Y,\rho)$ has the property forward parallel implies backward parallel, then the following are equivalent:
				
				\begin{enumerate}[(a)]
					\item $f$ is left K-Cauchy regular;
					\item for every real-valued forward uniformly continuous map $h:(Y,\rho)\rightarrow (\mathbb{R},|\cdot|)$,  we have $h\circ f:(X,d)\rightarrow (\mathbb{R},|\cdot|)$ is a left K-Cauchy regular map.
				\end{enumerate}
			\end{theorem}
			
			\begin{proof}
				$(a)\implies(b)$. It follows from the facts that every forward uniformly continuous function is left K-Cauchy regular and composition of left K-Cauchy regular maps is left K-Cauchy regular.\\
				$(b)\implies(a)$. On the contrary, suppose $f$ is not left K-Cauchy regular. So there exist sequences $(x_n)$ and $(y_n)$ in $(X,d)$ such that  $(x_n)\underset{R} {\overset{L}{\sim}}(y_n)$ but $(f(x_n))\underset{R} {\overset{L}{\nsim}}(f(y_n))$. Without loss of generality, let $(f(x_n)){\overset{L}{\nsim}}(f(y_n))$. Therefore, there exists $\epsilon >0$ such that for every $k \in \mathbb{N}$, there exist $n_k, m_k \in \mathbb{N}$ satisfying $n_k\geq m_k > k$ and $\rho (f(x_{m_{k}}),f(y_{n_{k}}))\geq\epsilon$. We can assume that $1 < m_1 \leq n_1 < m_2\leq n_2< \ldots$. So by Lemma \ref{l2}, we get an infinite subset $E$ of $\mathbb{N}$ such that $\rho^{s}(f(x_{m_{k}}),f(y_{n_{l}}))\geq \epsilon/4$ for every $k,l \in E$. Consequently, there exists an infinite subset $E_1$ of $E$ such that either $\rho(f(x_{m_{k}}),f(y_{n_{l}}))\geq \epsilon/4$ for all $k,l\in E_1$ or $\bar{\rho}(f(x_{m_{k}}),f(y_{n_{l}}))=\rho(f(y_{n_{l}}),f(x_{m_{k}}))\geq \epsilon/4$ for all $k,l\in E_1$. Take $A=\{f(x_{m_{k}}):k\in E_1\}$. If $\rho(f(x_{m_{k}}),f(y_{n_{l}}))\geq \epsilon/4$ for all $k,l\in E_1$, then define $h:(Y,\rho)\rightarrow (\mathbb{R},|\cdot|)$ as $h(x)=\rho (A,.)$; and if $\rho(f(y_{n_{l}}),f(x_{m_{k}}))\geq \epsilon/4$ for all $k,l\in E_1$, then define $h:(Y,\rho)\rightarrow (\mathbb{R},|\cdot|)$ as $h(x)=\rho (.,A)$. Now by Proposition \ref{funiformlycts}, in both the cases the map $h$ is forward uniformly continuous. However, $(h\circ f)(x_n)\underset{R} {\overset{L}{\nsim}}(h\circ f)(y_n)$. Consequently, $h\circ f$ is not left K-Cauchy regular. 
			\end{proof}

			We end this section with the following result that says a forward uniformly continuous function defined on a subset of a quasi-metric space can be extended to its forward closure. 
			\begin{theorem}
				Let $(X,d)$ and $(Y,\rho)$ be $T_1$ quasi-metric spaces with the property forward convergence implies backward convergence. If $(Y,\rho)$ is left K-Complete and $A\subseteq X$, then a forward uniformly continuous function $f:(A,d)\rightarrow (Y,\rho)$ can be uniquely extended forward uniformly to $Cl^+A$.
			\end{theorem}
			\begin{proof}
				Since every forward uniformly continuous function is left K-Cauchy regular, define $\bar{f}$ as in Theorem \ref{t3}. Now we show that $\bar{f}$ is forward uniformly continuous. Let $\epsilon >0$. Since $f$ is forward uniformly continuous on $A$, choose $\delta >0$ such that $\rho(f(a),f(b))<\epsilon/3$, whenever $a,b\in A$ and $d(a,b)<\delta$. Now if $x,y\in Cl^+A$ such that $d(x,y)<\delta$, then there exists $(x_n)$ and $(y_n)$ in $A$ such that $x_n\xrightarrow{f} x$ and $y_n\xrightarrow{f} y$. Since $(X,d)$ is $T_1$ quasi-metric spaces with the property forward convergence implies backward convergence and $x_n\xrightarrow{f} x$ and $y_n\xrightarrow{f} y$, there exists $n_0\in\mathbb{N}$ such that $d(x_n,y_n)-d(x,y)< \epsilon $ and $d(x,y)-d(x_n,y_n)< \epsilon $ for all $n\geq n_0$. Therefore, $\lim d(x_n,y_n)=d(x,y)$. As $\delta-d(x,y)>0$,  pick $N\in \mathbb{N}$ such that $d(x_n,y_n)<\delta$ for every $n\geq N$. Then $\rho(f(x_n),f(y_n))<\epsilon/3$  for every $n\geq N$. Consequently, there exists $n_\epsilon\in \mathbb{N}$ such that $\rho(\bar{f}(x),\bar{f}(y))\leq \rho(\bar{f}(x),f(x_n))+\rho(f(x_n),f(y_n))+\rho(f(y_n),\bar{f}(y))<\epsilon$ for all $n\geq n_\epsilon$. Hence $\bar{f}$ is forward uniformly continuous on $Cl^+A$.
				
				Uniqueness of $\bar{f}$ follows from the fact that every forward uniformly continuous function is ff-continuous and the Remark \ref{R1}.
			\end{proof}

			\bibliographystyle{amsplain}
			\bibliography{Reference_file}
		\end{document}